# Correction of a theorem in the Wiener index on fuzzy graphs

## Masoud Ghods[1,*], Zahra Rostami [2]

[1]Department of Mathematics, Semnan University, Semnan 35131-19111, Iran, mghods@semnan.ac.ir
[2]Department of Mathematics, Semnan University, Semnan 35131-19111, Iran, zahrarostami.98@semnan.ac.ir

* Correspondence: mghods@semnan.ac.ir; Tel.: (09122310288)

**Abstract**

In an article entitled The **Counterexample of a theorem in "Wiener index of a fuzzy graph and application to illegal immigration networks"**, we have shown a few examples of incorrect proof of a theorem in the Wiener index. Now, in this article, we will prove it by correcting the desired theorem. We have also used the same method as the author to prove it. There are other methods of proof.

**Keywords:** Fuzzy graph; Connectivity index; Wiener index; Fuzzy tree; Fuzzy cycle

## 1. Introduction and Preliminaries

Considering that this article is in order to correct a case in article "Wiener index of a fuzzy graph and application to illegal immigration networks", so we ask the reader to refer to the main article. This article is also published in [3].

Let $S$ be a set. A fuzzy graph is a pair $G:(\sigma,\mu)$, where $\sigma$ is a fuzzy subset of a set $S$ and $\mu$ is a fuzzy subset of $S\times S$ such that $\mu(u,v)\leq \sigma(u)\wedge\sigma(v)$. Note that $\mu$ is called a fuzzy relation on $\sigma$. We assume that $S$ is finite and nonempty, $\mu$ is reflexive and symmetric. The underlying graph of $G:(\sigma,\mu)$ is denoted by $G^*:(\sigma^*,\mu^*)$, where $\sigma^*=\{v\in S|\sigma(v)>0\}$ and $\mu^*=\{(u,v)\in S\times S\,|\mu(u,v)>0\}$. We denote an element $(x,y)$ of $\mu^*$ by $xy$ and call it an edge of $G$. Elements of $\sigma^*$ are called nodes or vertices of the fuzzy graph $G$. if the corresponding $G^*$ is trivial, we say that the fuzzy graph $G:(\sigma,\mu)$ is trivial. A sequence of distinct vertices $u_0,u_1,\ldots,u_n$ where $\mu(u_{i-1}u_i)>0, for\ i=1,2,\ldots n$ is called a path $P$ of length $n$. The membership value of a weakest edge is defined to be strength of the path $P$. A path $P$ is called a cycle if $u_0=u_n$. $CONN_G(x,y)$ is used to denote the strength of connectedness between two vertices $x$ and $y$ and is defined as the maximum of strengths of all paths between $x$ and $y$. if the strength of a path $P$ is equal to $CONN_G(x,y)$, then $P$ is called a strongest $x-y$ path. For each $x,y\in\sigma^*$, if $CONN_G(x,y)>0$, then $G:(\sigma,\mu)$ is called a connected fuzzy graph. An edge $xy$ of a fuzzy graph $G:(\sigma,\mu)$ is called $\alpha$-strong if $\mu(xy)>CONN_{G-xy}(x,y)$. $xy\in\mu^*$ is called $\beta$-strong if $\mu(xy)=CONN_{G-xy}(x,y)$ and it is called a $\delta$-edge if $\mu(xy)<CONN_{G-xy}(x,y)$. An edge is *strong* if it is either $\alpha$-strong or $\beta$-strong. A path $P$ is called a *strong path* if all of its edges are strong.

A connected fuzzy graph $G: (\sigma, \mu)$ is called a *fuzzy tree* if it has a spanning fuzzy subgraph $F = (\sigma, \upsilon)$ which is a tree, where for all $xy$ not in $F$, there exists a path from $x$ to $y$ in $F$, whose strength is more than $\mu(xy)$. Also $F$ is the unique maximum spanning tree ($MST$) of $G$.

**Definition 1.1**. [2] Let $G: (\sigma, \mu)$ be a Fuzzy graph. The *Connectivity index* $(CI)$ of G defined by

$$CI(G) = \sum_{u,v \in \sigma^*} \sigma(u)\sigma(v)CONN_G(u,v),$$

Where $CONN_G(u,v)$ is the strength of connectedness between $u$ and $v$.

**Definition 1.2**. [1] Let $G: (\sigma, \mu)$ be a Fuzzy graph. The *Wiener index* $(WI)$ of $G$ defined by

$$WI(G) = \sum_{u,v \in \sigma^*} \sigma(u)\sigma(v)d_s(u,v),$$

Where $d_s(u,v)$ is the minimum sum of weights of geodesics from $u$ to $v$.

**Definition 1.3**. [1] Let $G: (\sigma, \mu)$ be a fuzzy graph and $x, y \in \sigma^*$. A strong path $P$ from $x$ to $y$ is called a geodesic if there is no shorter strong path from $x$ to $y$. the weight of a geodesic is the sum of membership values of all edges in the geodesic.

**Theorem 1.4.** [2] Let $C_n$ be a fuzzy cycle. Then it is saturated if and only if the following two conditions are satisfied.

(i)     $n = 2k$, where k is an integer.
(ii)     $\alpha$-strong and $\beta$-strong edges appear alternatively in $C_n$.

**Theorem 1.5.** [2] Let $G: (\sigma, \mu)$ be a fuzzy tree and $F$ be the corresponding $MST$ of $G$. then $CI(F) = CI(G)$.

**Theorem 1.6.** [1] Let $G: (\sigma, \mu)$ be a fuzzy tree with $F$. Then $WI(G) = WI(F)$.

## 2. Results

Here the authors used Theorems 2.5 and 2.6 to achieve this equal, but as we have seen, this equal is not established because of the inequality of geodesic and connectivity in paths with length greater than one.

**Theorem* [1].** Let $G: (\sigma, \mu)$ be a saturated fuzzy cycle with $C^* = C_n$ of length $n$ such that each $\alpha$-strong edge has strength $\kappa$ and that of each $\beta$-strong edge is $\eta$, then

$$WI(G) = \frac{n[(n+3)^2 - 6]}{16}(\kappa + \eta).$$

**Proof.** Let $G: (\sigma, \mu)$ be a saturated fuzzy cycle. Then alternate edges of $C_n$ has strength $\kappa$ and $\eta$ respectively. By Theorem 1.4. $n$ is an even number. Note that each edge of a fuzzy cycle is strong. The maximum length of a geodesic in $C_n$ is $\frac{n}{2}$.

For $1 \leq k \leq \frac{n}{2}$, define

$$P_k = \{(u,v) \in \sigma^* \times \sigma^* \mid length\ of\ the\ geodesic\ between\ u\ and\ v\ is\ k\}.$$

There are $\frac{n}{2}$ pairs of vertices $(u,v)$ such that the length of the geodesic between $u$ and $v$ is $\frac{n}{2}$. $d_s(u,v) = \frac{n}{4}(\kappa + \eta)$ for any $(u,v) \in P_{\frac{n}{2}}$. Thus $\sum_{u,v \in P_{\frac{n}{2}}} d_s(u,v) = \frac{n}{4}(\kappa + \eta)\frac{n}{2} = \frac{n^2}{8}(\kappa + \eta)$.

Note that $|\sigma^*| = n$. Corresponding to each vertex $v \in \sigma^*$, there are precisely two vertices at a distance k from v. thus there are 2n pairs of vertices. If we avoid repetition, there will be n pairs $(x,y) \in \sigma^* \times \sigma^*$ such that length of the geodesic from x to y is k.

For $1 \leq k < \frac{n}{2}$, k even and $(u,v) \in P_k$, there are $\frac{k}{2}$ number of $\alpha$ and $\beta$-strong edges in the geodesic of length k from u to v. thus $\sum_{u,v \in P_k} d_s(u,v) = \frac{nk}{2}(\kappa + \eta)$. Let $1 \leq k < \frac{n}{2}$, k odd and $u \in \sigma^*$. Let v, w be the vertices at a distance k from u. one of the geodesics from u of length k contains $\frac{k+1}{2}$ $\alpha$-strong edges and $\frac{k-1}{2}$ $\beta$-strong edges. The other one contain $\frac{k+1}{2}, \frac{k-1}{2}$ of $\beta$ and $\alpha$-strong edges respectively.

$$\sum_{u,v \in P_k} d_s(u,v) = \frac{n}{2}\left[\frac{k+1}{2}\kappa + \frac{k-1}{2}\eta\right] + \frac{n}{2}\left[\frac{k+1}{2}\eta + \frac{k-1}{2}\kappa\right] = \frac{n}{2}[k\kappa + k\eta] = \frac{nk}{2}[\kappa + \eta].$$

Thus,

$$WI(G) = \frac{n^2}{8}(\kappa + \eta) + \sum_{1 \leq k < \frac{n}{2}} \frac{nk}{2}(\kappa + \eta) = (\kappa + \eta)\frac{n}{2}\left[\frac{n}{4} + \sum_{1 \leq k \leq \frac{n+1}{2}} k\right]$$

$$= \frac{n[(n+3)^2 - 6]}{16}(\kappa + \eta).$$

□

**Example 2.2**. Let $G: (\sigma, \mu)$ be a saturated fuzzy cycle with $C^* = C_4$ of length 4 such that each $\alpha$-strong edge has strength $\kappa$ and that of each $\beta$-strong edge is $\eta$, with $\sigma^* = \{a,b,c,d\}$, $\sigma(x) = 1$ for every $x \in \sigma^*$. Let $\mu(ab) = \mu(cd) = \kappa$ and $\mu(bc) = \mu(ad) = \eta$. Then

$$WI(G) = d_s(a,b) + d_s(a,c) + d_s(a,d) + d_s(b,c) + d_s(b,d) + d_s(c,d)$$
$$= \kappa + (\kappa + \eta) + \eta + \eta + (\eta + \kappa) + \kappa = 4\kappa + 4\eta = 4(\kappa + \eta).$$

But, by use *Theorem\**,

$$WI(G) = \frac{n[(n+3)^2 - 6]}{16}(\kappa + \eta) = \frac{4[(4+3)^2 - 6]}{16}(\kappa + \eta) = \frac{4[49 - 6]}{16}(\kappa + \eta)$$
$$= \frac{43}{4}(\kappa + \eta).$$

As can be seen, the value obtained is not correct.

Now, in the example blow let n=6.

**Example 2.3.** Let $G: (\sigma, \mu)$ be a saturated fuzzy cycle with $C^* = C_6$ of length 4 such that each $\alpha$-strong edge has strength $\kappa$ and that of each $\beta$-strong edge is $\eta$, with $\sigma^* = \{a, b, c, d, e, f\}$, $\sigma(x) = 1$ for every $x \in \sigma^*$. Let $\mu(ab) = \mu(cd) = \mu(ef) = \kappa$ and $\mu(bc) = \mu(de) = \mu(af) = \eta$.

|   | a | b | c | d | e | f |
|---|---|---|---|---|---|---|
| a | - | $\kappa$ | $\kappa + \eta$ | $2\eta + \kappa$ | $\eta + \kappa$ | $\eta$ |
| b | $\kappa$ | - | $\eta$ | $\eta + \kappa$ | $2\eta + \kappa$ | $\kappa + \eta$ |
| c | $\kappa + \eta$ | $\eta$ | - | $\kappa$ | $\kappa + \eta$ | $2\eta + \kappa$ |
| d | $2\eta + \kappa$ | $\eta + \kappa$ | $\kappa$ | - | $\eta$ | $\eta + \kappa$ |
| e | $\eta + \kappa$ | $2\eta + \kappa$ | $\kappa + \eta$ | $\eta$ | - | $\kappa$ |
| f | $\eta$ | $\kappa + \eta$ | $2\eta + \kappa$ | $\eta + \kappa$ | $\kappa$ | - |

Then

$$WI(G) = d_s(a,b) + d_s(a,c) + d_s(a,d) + d_s(a,e) + d_s(a,f) + d_s(b,c) + d_s(b,d) + d_s(b,e)$$
$$+ d_s(b,f) + d_s(c,d) + d_s(c,e) + d_s(c,f) + d_s(d,e) + d_s(d,f) + d_s(e,f)$$
$$= \kappa + (\kappa + \eta) + (2\eta + \kappa) + (\eta + \kappa) + (\eta) + (\eta) + (\eta + \kappa) + (2\eta + \kappa) + (\kappa + \eta)$$
$$+ (\kappa) + (\kappa + \eta) + (2\eta + \kappa) + (\eta) + (\eta + \kappa) + (\kappa) = 12\kappa + 15\eta.$$

But, by use *Theorem\**, we have

$$WI(G) = \frac{n[(n+3)^2 - 6]}{16}(\kappa + \eta) = \frac{6[(6+3)^2 - 6]}{16}(\kappa + \eta) = \frac{6[81 - 6]}{16}(\kappa + \eta)$$
$$= \frac{6 \times 75}{16}(\kappa + \eta)$$

As can be seen, the value obtained is not correct.

Now, in the example blow let n=8.

**Example 2.4.** Let $G: (\sigma, \mu)$ be a saturated fuzzy cycle with $C^* = C_8$ of length 4 such that each $\alpha$-strong edge has strength $\kappa$ and that of each $\beta$-strong edge is $\eta$, with $\sigma^* = \{a, b, c, d, e, f, g, h\}$, $\sigma(x) = 1$ for every $x \in \sigma^*$. Let $\mu(ab) = \mu(cd) = \mu(ef) = \mu(gh) = \kappa$ and $\mu(bc) = \mu(de) = \mu(fg) = \mu(ha) = \eta$.

|   | a | b | c | d | e | f | g | h |
|---|---|---|---|---|---|---|---|---|
| a | - | $\kappa$ | $\kappa + \eta$ | $\eta + 2\kappa$ | $2\eta + 2\kappa$ | $2\eta + \kappa$ | $\kappa + \eta$ | $\eta$ |
| b | $\kappa$ | - | $\eta$ | $\eta + \kappa$ | $2\eta + \kappa$ | $2\kappa + 2\eta$ | $2\kappa + \eta$ | $\kappa + \eta$ |
| c | $\kappa + \eta$ | $\eta$ | - | $\kappa$ | $\kappa + \eta$ | $\eta + 2\kappa$ | $2\kappa + 2\eta$ | $2\eta + \kappa$ |
| d | $\eta + 2\kappa$ | $\eta + \kappa$ | $\kappa$ | - | $\eta$ | $\eta + \kappa$ | $2\eta + \kappa$ | $2\kappa + 2\eta$ |
| e | $2\eta + 2\kappa$ | $2\eta + \kappa$ | $\kappa + \eta$ | $\eta$ | - | $\kappa$ | $\kappa + \eta$ | $2\kappa + \eta$ |
| f | $2\eta + \kappa$ | $2\kappa + 2\eta$ | $\eta + 2\kappa$ | $\eta + \kappa$ | $\kappa$ | - | $\eta$ | $\eta + \kappa$ |
| g | $\kappa + \eta$ | $2\kappa + \eta$ | $2\kappa + 2\eta$ | $2\eta + \kappa$ | $\kappa + \eta$ | $\eta$ | - | $\kappa$ |
| h | $\eta$ | $\kappa + \eta$ | $2\eta + \kappa$ | $2\kappa + 2\eta$ | $2\kappa + \eta$ | $\eta + \kappa$ | $\kappa$ | - |

Then
$$WI(G) = 32\kappa + 32\eta = 32(\kappa + \eta).$$

But, by use *Theorem\**, we have
$$WI(G) = \frac{n[(n+3)^2 - 6]}{16}(\kappa + \eta) = \frac{8[(8+3)^2 - 6]}{16}(\kappa + \eta) = \frac{8[121 - 6]}{16}(\kappa + \eta)$$
$$= \frac{8 \times 115}{16}(\kappa + \eta).$$

As can be seen, the value obtained is not correct.
We will now modify the theorem\* as follows.

**Theorem 2.1.** Let $G: (\sigma, \mu)$ be a saturated fuzzy cycle with $C^* = C_n$ of length $n$ such that each $\alpha$-strong edge has strength $\kappa$ and that of each $\beta$-strong edge is $\eta$, then
$$\begin{cases} \frac{n^3}{16}(\kappa + \eta) & n = 4m \ (m = 1, 2, \dots) \\ \frac{n(n^2 - 4)}{16}\kappa + \frac{n(n^2 + 4)}{16}\eta & O.W \end{cases}$$

Before we prove theorem, let's look at the examples above.
For $n = 4$, we have

$$\frac{n^3}{16}(\kappa + \eta) = \frac{4^3}{16}(\kappa + \eta) = 4(\kappa + \eta).$$

For $n = 6$,

$$\frac{n(n^2-4)}{16}\kappa + \frac{n(n^2+4)}{16}\eta = \frac{6(6^2-4)}{16}\kappa + \frac{6(6^2+4)}{16}\eta = \frac{6 \times 32}{16}\kappa + \frac{6 \times 40}{16}\eta$$
$$= 12\kappa + 15\eta.$$

For $n = 8$,

$$\frac{n^3}{16}(\kappa + \eta) = \frac{8^3}{16}(\kappa + \eta) = 32(\kappa + \eta).$$

We use the author's method in the main article to prove it.

**Proof theorem 2.1.** Let $G: (\sigma, \mu)$ be a saturated fuzzy cycle. Then alternate edges of $C_n$ has strength $\kappa$ and $\eta$ respectively. By Theorem 1.4. $n$ is an even number. Note that each edge of a fuzzy cycle is strong. The maximum length of a geodesic in $C_n$ is $\frac{n}{2}$.

For $1 \leq k \leq \frac{n}{2}$, define
$$P_k = \{(u, v) \in \sigma^* \times \sigma^* \mid length\ of\ the\ geodesic\ between\ u\ and\ v\ is\ k\}.$$

There are $\frac{n}{2}$ pairs of vertices $(u, v)$ such that the length of the geodesic between $u$ and $v$ is $\frac{n}{2}$. If $n = 4m$, for $m = 1, 2, 3, ....$ ($n$ is a factor of 4) We have $d_s(u, v) = \frac{n}{4}(\kappa + \eta)$ for any $(u, v) \in P_{\frac{n}{2}}$, And $\frac{(n-2)}{4}\kappa + \frac{(n+2)}{4}\eta$, for elsewhere($n$ is not a factor of 4) with any $(u, v) \in P_{\frac{n}{2}}$. Thus $\sum_{u,v \in P_{\frac{n}{2}}} d_s(u, v) = \frac{n}{4}(\kappa + \eta)\frac{n}{2} = \frac{n^2}{8}(\kappa + \eta)$ ( $n$ is a factor of 4), and $\frac{n}{2}(\frac{(n-2)}{4}\kappa + \frac{(n+2)}{4}\eta)$. Note that $|\sigma^*| = n$. Corresponding to each vertex $v \in \sigma^*$, there are precisely two vertices at a distance $k$ from $v$. thus there are $2n$ pairs of vertices. If we avoid repetition, there will be n pairs $(x, y) \in \sigma^* \times \sigma^*$ such that length of the geodesic from $x$ to $y$ is $k$.

For $1 \leq k < \frac{n}{2}$, k even and $(u, v) \in P_k$, there are $\frac{k}{2}$ number of $\alpha$ and $\beta$-strong edges in the geodesic of length $k$ from $u$ to $v$. thus $\sum_{u,v \in P_k} d_s(u, v) = \frac{nk}{2}(\kappa + \eta)$. Let $1 \leq k < \frac{n}{2}$, k odd and $u \in \sigma^*$. Let $v, w$ be the vertices at a distance $k$ from $u$. One of the geodesics from $u$ of length $k$ contains $\frac{k+1}{2}$ $\alpha$-strong edges and $\frac{k-1}{2}$ $\beta$-strong edges. The other one contain $\frac{k+1}{2}, \frac{k-1}{2}$ of $\beta$ and $\alpha$-strong edges respectively. Then

$$\sum_{u,v \in P_k} d_s(u,v) = \frac{n}{2}\left[\frac{k+1}{2}\kappa + \frac{k-1}{2}\eta\right] + \frac{n}{2}\left[\frac{k+1}{2}\eta + \frac{k-1}{2}\kappa\right] = \frac{n}{2}[k\kappa + k\eta] = \frac{nk}{2}[\kappa + \eta].$$

Thus, For $n = 4m, m = 1, 2, 3, \ldots$

$$WI(G) = \frac{n^2}{8}(\kappa + \eta) + \sum_{1 \le k < \frac{n}{2}} \frac{nk}{2}(\kappa + \eta) = (\kappa + \eta)\frac{n}{2}\left[\frac{n}{4} + \sum_{1 \le k \le \frac{n}{2}-1} k\right]$$

$$= \frac{n}{2}\left[\frac{n}{4} + (\frac{n}{2} - 1) \times \frac{1 + \frac{n}{2} - 1}{2}\right](\kappa + \eta) = \frac{n}{2}\left[\frac{n}{4} + \frac{n-2}{2} \times \frac{n}{4}\right](\kappa + \eta)$$

$$= \frac{n}{2}\left[\frac{n}{4} + \frac{n^2 - 2n}{8}\right](\kappa + \eta) = \frac{n}{2}\left[\frac{n^2}{8}\right](\kappa + \eta) = \frac{n^3}{16}(\kappa + \eta).$$

For $n$ where $n$ is not a factor of $4$

$$WI(G) = \frac{n}{2}\left(\frac{(n-2)}{4}\kappa + \frac{(n+2)}{4}\eta\right) + \sum_{1 \le k < \frac{n}{2}} \frac{nk}{2}(\kappa + \eta)$$

$$= \frac{n}{2}\left(\frac{(n-2)}{4}\kappa + \frac{(n+2)}{4}\eta\right) + (\kappa + \eta)\frac{n}{2}\left[\sum_{1 \le k \le \frac{n}{2}-1} k\right]$$

$$= \frac{n}{2}\left(\frac{(n-2)}{4}\kappa + \frac{(n+2)}{4}\eta\right) + \frac{n}{2}\left[(\frac{n}{2} - 1) \times \frac{1 + \frac{n}{2} - 1}{2}\right](\kappa + \eta)$$

$$= \frac{n}{2}\left(\frac{(n-2)}{4}\kappa + \frac{(n+2)}{4}\eta\right) + \frac{n}{2}\left[\frac{n(n-2)}{8}\right](\kappa + \eta)$$

$$= \frac{n(n^2 - 4)}{16}\kappa + \frac{n(n^2 + 4)}{16}\eta.$$

□

**Refrences**